\documentclass[preprint]{elsarticle}
\usepackage{graphicx}
\usepackage{subcaption}
\usepackage{amsmath,amsthm,amsfonts,mathtools}
\usepackage{derivative}
\usepackage{algorithm}
\usepackage{algcompatible}
\usepackage{algpseudocode}
\DeclareMathOperator{\diag}{diag}

\DeclareMathOperator{\span1}{span}

\usepackage{listings}
\usepackage{subdepth} 
\usepackage{mathtools}
\newtheorem{theorem}{Theorem}
\newtheorem{proposition}[theorem]{Proposition}

\newcommand{\B}[1]{\mbox{\boldmath $#1$}}

\newcommand{\iu}{{i\mkern1mu}}

\begin{document}

\begin{frontmatter}

\title{Efficient Inversion of  Matrix $\phi$-Functions of Low Order}
\author[1]{Luca Gemignani\fnref{fn1}}
\ead{luca.gemignani@unipi.it}
\address{Dipartimento di Informatica, Largo B. Pontecorvo, 3, Pisa,
  56127, Italy}
\fntext[fn1]{The author is partially supported by INDAM/GNCS and by the project PRA\_2020\_61 of the University of Pisa.}

\begin{abstract}
  The paper is concerned with efficient numerical methods for solving a linear system $\phi(A) \B x=\B b$, where  $\phi(z)$ is a $\phi$-function and $A\in \mathbb R^{N\times N}$.
  In particular  in this work we are interested in  the computation of  ${\phi(A)}^{-1}\B b$  for  the case where $\phi(z)=\phi_1(z)=\displaystyle\frac{e^z-1}{z}, \quad \phi(z)=\phi_2(z)=\displaystyle\frac{e^z-1-z}{z^2}$.
  Under suitable conditions on the spectrum of $A$   we design  fast algorithms    for computing   both ${\phi_\ell(A)}^{-1}$  and ${\phi_\ell(A)}^{-1}\B b$  
  based on  Newton's iteration  and  Krylov-type methods, respectively. 
  Adaptations  of these  schemes for structured matrices
  are considered. In particular the cases of banded and more generally quasiseparable matrices are investigated.
  Numerical results are presented to show the effectiveness of our proposed algorithms. 
\end{abstract}


\end{frontmatter}

\section{Introduction}
Efficient numerical methods for computing the action of matrix $\phi$-functions are  of growing interest for the application of
exponential integrators  in the solution of stiff systems of differential equations (compare \cite{Minchev2005ARO,HL,HO,NW,GRT,CCZ} and the references given therein).
The computation of the inverse of  matrix $\phi$-functions or, equivalently,  the design of fast linear solvers for  matrix $\phi$-functions is useful in the solution of
related inverse problems. 

A fast efficient numerical method  for computing $\psi_1(A)$ and $\psi_1(A)\B b$  with  $\psi_1(z)=1/\phi_1(z)$, $\phi_1(z)=\displaystyle\frac{e^z-1}{z}$, $A\in \mathbb R^{N\times N}$,   has been presented in \cite{BEG1, BEG2}.
The method exploits a partial
fraction decomposition of the
meromorphic function  $\psi_1(z)$ and it is  particularly suited for the application to structured matrices  for which fast linear solvers exist.  The same  approach cannot be extended to other  functions
$\psi_\ell(z)=1/\phi_\ell(z)$  with $\ell>1$  due to the lack of explicit closed--form expressions of their poles. A numerical investigation of the poles of $\psi_2(A)$   is carried out in  \cite{KT}.

The computation of 
 $\psi_2(A)$ and $\psi_2(A)\B b$  with  $\psi_2(z)=1/\phi_2(z)$, $\phi_2(z)=\displaystyle\frac{e^z-1 -z}{z^2}$, $A\in \mathbb R^{N\times N}$, is also of relevant interest.  We describe hereafter two  applications. 
\begin{enumerate}
    \item {\tt A nonlocal inverse problem.} Consider the  nonlocal inverse  problem defined as follows: We seek the vector $\B g\in \mathbb R^N$ and the function $\B u=\B u(t)\colon [0,T]\rightarrow \mathbb R^N$ such that 
    \[
    \B u'(t)=A\B u(t) +\B g,  \quad \B u(0)=\B u_0, \ \  
    \displaystyle\frac{1}{T}\int_0^T\B u(t)dt=\B u_1
    \]
    with  $A\in \mathbb R^{N\times N}$, $\B u_0, \B u_1\in \mathbb R^N$.  The existence and uniqueness  of the solution in an abstract Banach space is   investigated in \cite{KT} by relying upon the theory developed in \cite{TEY}.  Under quite general assumptions it is shown that $\B u(t)$ can be taken  of the form 
    $\B u(t)= e^{tA} \B u_0 +\int_0^t e^{(t-s)A} \B g  ds$. By imposing the integral condition we find that the vector $\B g$  solves 
    \begin{equation}\label{init1}
    \B u_1=\phi_1(TA)\B u_0 +T \phi_2(TA)\B g.
    \end{equation}
Sufficient  conditions  for the solvability of the system can  be expressed in terms of the eigenvalues of $A$ \cite{TEY}. 
\item {\tt Two-point inverse problems.} The computation of the unknown  parameter $\B p$  in the local boundary  value problem   \cite{POV,Suh} 
\begin{equation}\label{suh1}
\B u'(t)=A \B u(t) +t \B p, \quad \B u(0)=\B q,  \ \ \B u(1)=\B g,
\end{equation}
$A\in \mathbb R^{N\times N}$, $t\in [0,1]$, amounts  to solve a linear system with a matrix  $\phi$-function as coefficient matrix. 
Exponential-type methods  \cite{HO}  can be used for  integrating  \eqref{suh1}. These methods follow from  from the   exact formula
\begin{equation}\label{exint00}
\B u(t + h) =  e^{hA} \B u(t) +e^{hA}\int_0^h e^{-As}
\B p (t+s)ds
\end{equation}
which  gives 
\[
\B u(1)= e^{A} \B u(0) +e^{A}\int_0^{1} e^{-As}\B p \ s  \ ds.
\]
By  computing the integral we  find
\begin{equation}\label{invex}
\B u(1) = \B g=  e^{A}\B q + \phi_2(A)\B p.
\end{equation}
This  relation makes  possible to compute $\B p$   by solving a linear system with coefficient matrix  $M=\phi_2(A)$.  
\end{enumerate}

In this contribution we address the computation of
$\psi_{\ell+1}(A)$   and $\psi_{\ell+1}(A)\B b$ with $\psi_{\ell}(z)=1/\phi_\ell(z)$ and $\ell>0$.
Our extension  relies on the Newton iteration for computing the inverse of a matrix.   For  review of this method  see \cite{PR,Pan}.  This tool has been successfully applied in  \cite{BM,PRW} for the inversion of matrices having a displacement rank structure.  
To make Newton's method work for matrix inversion, an initial approximate inverse $X_0$ of $B=\phi_{\ell+1}(A)$ is required.  Then, it is easily seen that  the   intermediate approximations of $\psi_{\ell+1}(A)$ generated  by Newton's iteration can be  expressed as a polynomial of  $X_0 B $ thus providing the link to the development of Krylov-type methods  for computing $\psi_{\ell+1}(A)\B b$. 

More specifically,   we  first identify regions $\Omega \subset \mathbb  C$  of the complex plane such  that $\left|1-\displaystyle\frac{\phi_{\ell+1}(z)}{ \phi_{\ell}(z)}\right|<1$, $\ell\geq 1$,  for all $z\in \Omega$.
Then, we show that  if the  eigenvalues of $A$ lie in $\Omega$,    the Newton iteration  applied for the inversion
  of $B=\phi_{\ell+1}(A)$, $\ell\geq 1$, with starting point $X_0=\psi_\ell(A)=(\phi_\ell(A))^{-1}$  is
  quadratically converging to the inverse matrix of $B$.  This means that given a method to compute  $X_0=\psi_\ell(A)$ we can  apply the Newton iteration for approximating
  $B_0^{-1}=\psi_{\ell+1}(A)$. 
  Moreover, since  the iterative scheme only requires  matrix multiplications it is amendable for fast implementations using
structured representations of the matrices involved.  In particular, fast adaptations  for  both displacement structured and quasiseparable matrices  can be devised.   Approximate  compression techniques in the style of
\cite{BM} can also  be incorporated to take under control the growth of  displacement or  quasiseparable ranks.

The convergence results  for the  Newton   iteration  can also be exploited in a different  perspective.  It is shown that the approximation $X_k$ of $\psi_{\ell+1}(A)$  obtained after $k$ iterations  satisfies 
\[
X_k= p_k(\psi_\ell(A) \phi_{\ell+1}(A))\psi_\ell(A), \quad k>0, 
\]
where $p_k(z)$ is a polynomial  of degree $2^k-1$.  If the convergence is very rapid, then  the solution of the  linear system $\phi_{\ell+1}(A) \B x=\B b$ can be approximated  efficiently by means of a Krylov-type method
like GMRES applied  for solving the 
equivalent system $\psi_\ell(A)\phi_{\ell+1}(A) \B x=\psi_\ell(A)\B b$. When  the  eigenvalues of $A$ lie in $\Omega$ then the  convergence of GMRES  applied to this system  follows  from  the results in \cite{saadb}
(see Proposition 6.32  and its generalizations).
The  paramount advantage of  such  a Krylov-based  approach is that matrix-by-vector multiplications are only 
required  to  find an approximation $X_k \B b $ of the vector $\B x=\psi_{\ell+1}(A)\B b$.  In particular, for $\ell=1$ the projection method  only involves   products of the form  $\psi_1(A) \B v$ which can be computed   using the methods introduced in \cite{BEG1, BEG2}.

In principle,  the proposed  schemes  can be applied recursively for evaluating  $\psi_{\ell+1}(A)$   or  $\psi_{\ell+1}(A)\B b$, $\ell>0$, provided that  a method for computing $\psi_{1}(A)$   or  $\psi_{1}(A)\B b$ is available. Despite the generality of the approach, however,  based on application  and  numerical issues in this work  we focus on the  case $\ell\in\{0,1\}$, or at least $\ell$ small in value. 

The  paper is organized as follows. In Section \ref{two} we  recall some   preliminaries   on both  $\phi-$ and $\psi-$functions,   the Newton  iteration for matrix inversion and  its connection with Krylov-type methods. In Section \ref{three}    we
analyze theoretical and computational properties of  Newton's iteration for the inversion of matrix 
$\phi$-functions. In Section \ref{threenew}   we devise   a Krylov-type method  for computing  the action of these inverses on a vector.   In Section \ref{five} we present the results of numerical experiments
illustrating the properties  of this method  whereas  conclusions  and future work are drawn in Section \ref{six}.

\section{Preliminaries}\label{two}
The $\phi$-functions are  entire functions defined for scalar arguments by the integral representation
\begin{equation}\label{intrep}
\phi_0(z)=e^{\displaystyle{z}}, \quad \phi_\ell(z)=\frac{1}{(l-1)!}\int_0^1 e^{\displaystyle{(1-\theta)z}} \theta ^{\ell-1} d\theta, \quad \ell\geq 1, \ z\in \mathbb C.
\end{equation}
The  $\phi$-functions satisfy the recurrence relation
\begin{equation}\label{recurrence}
\phi_{\ell}(z)=z \phi_{\ell+1}(z) +\frac{1}{l!}, \quad \ell\geq 0,
\end{equation}
and have the Taylor expansion
\[
\phi_{\ell}(z)=\sum_{k=0}^\infty \frac{z^k}{(k+\ell)!}, \quad  \ell\geq 0.
\]
This latter  can be extended to a matrix argument by setting for any $A\in \mathbb C^{N\times N}$
\[
\phi_{\ell}(A)=\sum_{k=0}^\infty \frac{A^k}{(k+\ell)!}, \quad  \ell\geq 0.
\]
The function  $\psi_\ell(z)$, $\ell\geq 0$,  is a meromorphic
function defined as the reciprocal of $\phi_\ell(z)$, that is,
\[
\psi_\ell(z)=\phi_\ell(z)^{-1}, \quad \ell\geq 0.
\]
Explicit series expansions are only known for $\psi_1(z)$.  It holds [\cite{AS}, formula 23.1.1]
  \begin{equation}\label{pws}
   \psi_1(z)={\phi_1(z)}^{-1}=\sum_{k=0}^{+\infty} \frac{B_k}{k!}z^k, \quad |z|< 2 \pi, 
    \end{equation}
    where $B_k$ denotes  the $k$th  Bernoulli number.  A different rational representation is derived  in \cite{BEG2}.  For any fixed $n>0$  we have 
\begin{equation}\label{ratapps}
 \psi_1( z) =f_n(z) +
 2(-1)^{n}\displaystyle\sum_{k=1}^{\infty}\left(\frac{z}{2 \pi}\right)^{2(n+1)}
    \frac{1}{k^{2n}}(\left(\frac{z}{2 \pi}\right)^2+ k^2)^{-1},
\end{equation}
    where 
 \[
  f_n(z)=  1 -\frac{1}{2} z+ 
\displaystyle\sum_{i=0}^{n-1} z^{2(i+1)} \frac{B_{2(i+1)}}{(2 (i+1))!}.
\]
The series  on the rhs of \eqref{ratapps} converges uniformly to $\psi_1( z)$  over any compact set $\mathcal K\subset \mathbb C\setminus \pm  2 \pi\iu \mathbb N$.  The polynomial 
contribution $f_n(z)$ is a partial sum of the power series expansion  \eqref{pws} aimed to improve the accuracy of the approximation around the  removable singularity at the origin  in the complex plane.  Relation \eqref{ratapps} provides a family of mixed polynomial/rational approximations of $\psi_1(A)$ of the form 
\begin{equation}\label{ratappsnm}
\psi_1(A)\simeq r_{n,m}(A)=f_n(A) + 2(-1)^{n}\displaystyle\sum_{k=1}^{m}\left(\frac{A}{2 \pi}\right)^{2(n+1)}
    \frac{1}{k^{2n}}(\left(\frac{A}{2 \pi}\right)^2+ k^2)^{-1}.
\end{equation}

Newton's iteration \cite{PR,Pan} for the inversion of a nonsingular matrix $B\in \mathbb C^{N\times N}$  is defined by :
\begin{equation}\label{Newton}
X_0\in \mathbb C^{N\times N}, \quad X_{k+1}=2 X_k  - X_k  B   X_k , \ k\geq 0.
\end{equation}
From 
\[
I-X_{k+1} B= (I-X_kB)^2=(I-X_0B)^{2^k},
\]
we obtain that  Newton's iteration \eqref{Newton}  quadratically converges  to $B^{-1}$ provided that all eigenvalues of $R=I-X_0B$ have modulus less than 1. 

Observe that 
\[
X_1=2 X_0-X_0  B   X_0= (2 I -X_0 B) X_0=p_1(X_0B) X_0, 
\]
with $p_1(z)$  a polynomial of degree 1.  Inductively, we find that for $k>0$
\begin{equation}\label{commute}
X_k=2 p_{k-1}(X_0B)X_0  -  p_{k-1}(X_0B) X_0 B p_{k-1}(X_0B) X_0 = p_k(X_0B) X_0, 
\end{equation}
for a suitable polynomial $p_k(z)$ of degree $2^k-1$.
This means that the approximation $\B x_k=X_k \B b$  of the solution $\B x$
 of $B\B x=\B b$  satisfies 
 \[
 \B x_k=p_k(X_0 B)  X_0\B b, \quad k>0.
 \]
  Hence  $\B x_k$ belongs to the  $m-th$  Krylov subspace, $m=2^k-1$, 
  \[
  \mathcal K_m=\span1\{X_0\B b, (X_0 B)X_0\B b, \ldots, (X_0 B)^m X_0\B b\}. \]
  It  follows that  a Krylov-type method  might be used for solving   the equivalent system 
  \begin{equation}\label{add1}
   X_0 B \B x= X_0\B b.
   \end{equation}
   For instance, GMRES \cite{saadb} after $m$ iterations returns an approximation $\B v_k$ such that 
   \[
  \parallel  X_0\B b-X_0 B \B v_k\parallel_2 \leq \parallel  X_0\B b-X_0 B \B x_k\parallel_2 \leq \parallel X_0\parallel_2 \parallel  \B b-B \B x_k\parallel_2.
  \]
  A  precise convergence estimate for GMRES  applied for the solution of  \eqref{add1}  when  all   the eigenvalues of $R=I-X_0B$ have modulus less than 1   will be given in Proposition \ref{gmresp}.

\section{Newton  Iteration for the Inversion of Matrix $\phi$-Functions}\label{three}

In this section we  design an  iterative method  based on Newton's iteration for the inversion of matrix $\phi$-functions $\phi_\ell(A)$, $A\in \mathbb C^{N\times N}$, $\ell>1$.

Let  us first  suppose that $A\in \mathbb C^{N\times N}$   has real eigenvalues only, that is, $\lambda \in {\tt spec}(A) \Rightarrow  \lambda \in  \Omega=\mathbb R$. 
Observe that for real  arguments  ($z\in \mathbb R$)
from the integral representation \eqref{intrep}  it follows that  $\phi_\ell (z)>0$  and,  moreover,   $\phi_\ell (z)>  \phi_{\ell+1}(z)$ $\forall z\in \mathbb R$,  $\ell\geq 1$.  This means that
\[
0< 1-\frac{\phi_{\ell+1}(\lambda_i)}{ \phi_{\ell}(\lambda_i)}=          1-\phi_{\ell+1}(\lambda_i) \psi_\ell(\lambda_i)<1, \quad \forall \lambda_i \in {\tt spec}(A).
\]
The next result immediately follows.

\begin{proposition}\label{prop1}
  Let $A\in \mathbb C^{N\times N}$   be a matrix with   all  real eigenvalues.  Then for any $\ell\geq 0$ $\phi_\ell(A)$ is invertible. Moreover, the Newton iteration \eqref{Newton} applied for the inversion
  of $B=\phi_{\ell+1}(A)$, $\ell\geq 1$, with starting point $X_0=\psi_\ell(A)=(\phi_\ell(A))^{-1}$  is
  quadratically converging to the inverse matrix of $B$.
  \end{proposition}

The extension  of this result for matrices with  possibly complex eigenvalues requires some additional constraints.
If $z=a+\iu b$, $a,b\in \mathbb R$, $\iu^2=-1$,  is a complex number then from the integral representation  \eqref{intrep}  it is found that for $\ell\geq 1$
\[
\phi_\ell (z)=\frac{\displaystyle\int_0^1 e^{\displaystyle{\tau a}} \cos(\tau b) (1-\tau)^{\ell-1} d\tau + \iu \displaystyle\int_0^1 e^{\displaystyle{\tau a}} \sin(\tau b) (1-\tau)^{\ell-1} d\tau}{(l-1)!}.
\]
Under the auxiliary assumption  $b\in [-\pi/2,\pi/2]$ this implies that $\Re(\phi_\ell (z))>0$ and, hence, $\phi_\ell (z)\neq 0$.  In addition, the residual $r(z)=\phi_{\ell}(z)-\phi_{\ell+1}(z)$  also satisfies
\[
r(z)=\frac{\displaystyle\int_0^1 e^{\displaystyle{\tau a}} \cos(\tau b)\frac{\ell-1+\tau}{\ell} (1-\tau)^{\ell-1} d\tau + \iu \displaystyle\int_0^1 e^{\displaystyle{\tau a}} \sin(\tau b) \frac{\ell-1+\tau}{\ell}(1-\tau)^{\ell-1} d\tau}{(l-1)!}.
\]
It follows that
\[
|\Re(r(z))|<|\Re(\phi_\ell (z))|, \quad |\Im(r(z))|<|\Im(\phi_\ell (z))|
\]
and therefore
\[
\left|\frac{r(z)}{\phi_\ell (z)}\right|=\left|1-\frac{\phi_{\ell+1}(z)}{ \phi_{\ell}(z)}\right|<1.
\]
To sum up  we arrive at the following  extension of Proposition \ref{prop1}.

\begin{proposition}\label{prop2}
  Let $A\in \mathbb C^{N\times N}$   be a matrix with   all   eigenvalues lying in the strip $ \Omega=\mathbb R\times \iu [-\pi/2,\pi/2]$ in the complex plane.  Then,  for any $\ell\geq 0$,   $\phi_\ell(A)$ is invertible.
  Moreover, the Newton iteration \eqref{Newton} applied for the inversion
  of $B=\phi_{\ell+1}(A)$, $\ell\geq 1$, with starting point $X_0=\psi_\ell(A)=(\phi_\ell(A))^{-1}$  is
  quadratically converging to the inverse matrix of $B$.
  \end{proposition}

Differently from the case of real spectrum   some restrictions on the localization of the eigenvalues are  needed for general matrices.  Let us consider the
tridiagonal Toeplitz matrix  $T$ of order $N$ having subdiagonal, diagonal and superdiagonal entries given by $0.5$, $0$ and $-0.5$,  respectively.  The matrix has eigenvalues located on the imaginary axis in the interval
$\iu [-1,1]$.  In Table \ref{t1},   we report  the  computed spectral radius of $R=I-(\phi_1(A))^{-1}(\phi_2(A)$, where  $A= h^{-1}  N^2 T$  and $N=128$, for different values of $h$. 

\noindent
\begin{table}
		\begin{center}
			\begin{tabular}{|l|l|l|l|l|} \hline
			
			$h$ & $1$ &$ N$ & $N^2$ & $N^4$ \\ 	\hline
			$\rho(R)$ & 1.6852e+03 & 57.5590 & 0.5071& 0.5000  \\ \hline
			\end{tabular}
\end{center}
       \caption{ Values of the spectral radius $\rho(R)$ of $R=I-(\phi_1(A))^{-1}\phi_2(A)$ with $A=h^{-1} N^2 T$, $N=128$,  and  $T={\tt {gallery('tridiag', N, 0.5, 0, -0.5)}}$  for different values of $h=1, N, N^2, N^4$. \label{t1}} 
\end{table}

Computational interest in Newton's method   is especially due to the development of high-performance computing environments. The iterative scheme \eqref{Newton}  basically  require BLAS Level 3 routines which are  easily implemented in parallel on a parallel computing system  \cite{PVP,Veneva2019}. Moreover, it is  especially suited to  take advantage of the sparsity  and  the structural properties of the matrices involved.  The case of matrices having a displacement structure  has been considered in
\cite{BM,Pan,PRW}.   In the next subsection, we focus on the application of Newton's method  for  inverting  banded and more generally quasiseparable-type  matrices arising from the discretization of partial differential equations.

\subsection{Fast Adaptations for Structured Matrices}\label{four}

\begin{figure}
  \centering
  \subfloat[]{\includegraphics[width=0.5\textwidth]{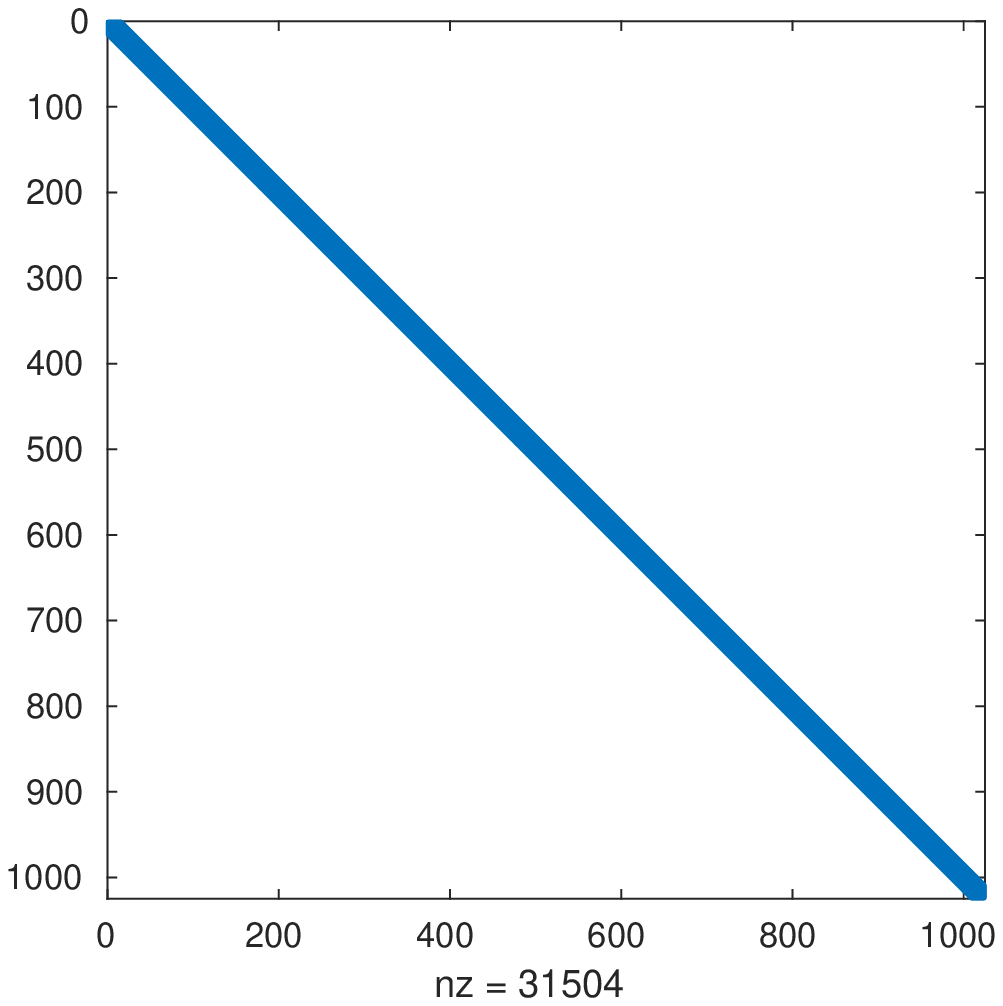}}
  \hfill
  \subfloat[]{\includegraphics[width=0.5\textwidth]{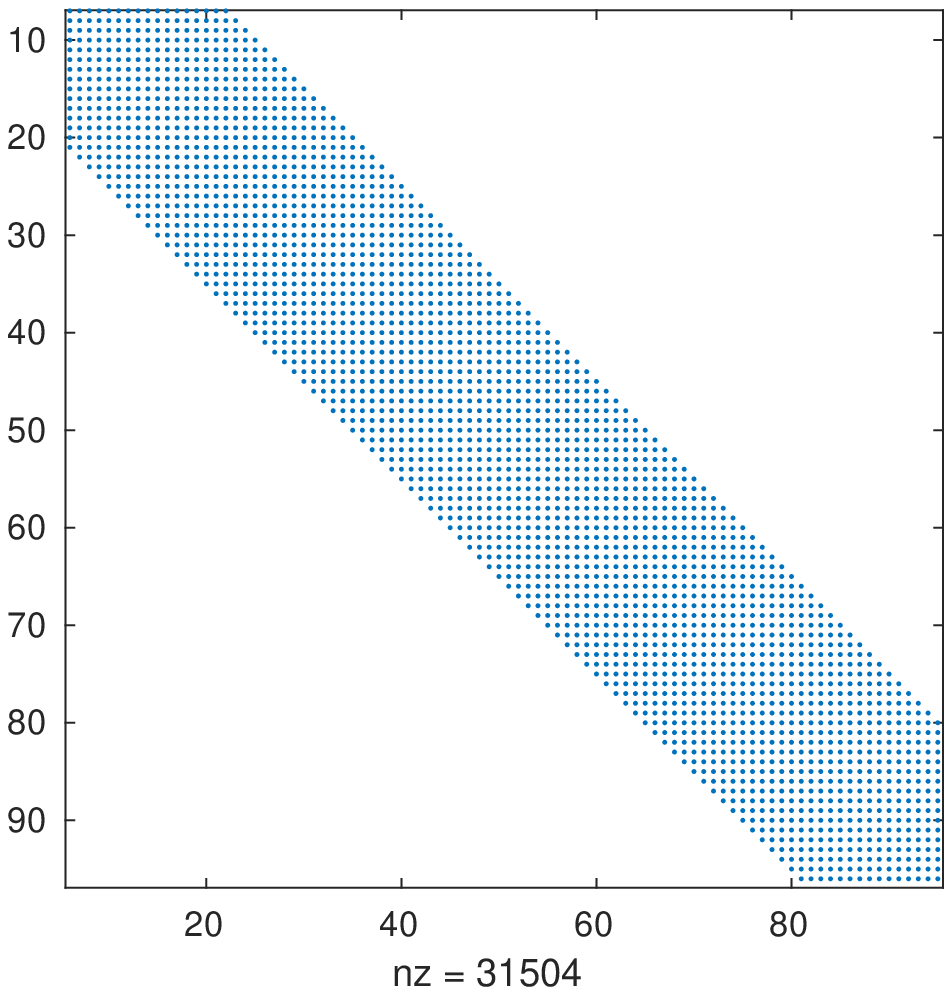}}
  \caption{"Spy" plots of $\psi_2(A)$ and its leading principal submatrix of order $95$ for $A$ being the 1D Laplacian matrix of order $N=1024$.}
\label{f3}
\end{figure}

We begin by observing that  for a given banded matrix $A$   
 the matrix $B=\phi_\ell(A)$  or  $B=\psi_\ell(A)$, $\ell\geq 1$, generally inherits the banded structure of  $A$ in some  approximate  way.  For instance, in Figure \ref{f3}  we illustrate the  "spy" plots of   $\psi_2(A)$ and its leading principal submatrix of order $95$ when  $A$  is the 1D Laplacian matrix of order $N=1024$. The  threshold  value is set to  $1.0e-14$.  The  exact tridiagonal structure of $A$  results into an approximate banded structure of $B$.  Precise mathematical statements depend on quantities that are  hard  to compute  and  typically yield very pessimistic estimates (compare with \cite{HLARGE} for the matrix exponential, the review \cite{Benzi} for  more  general analytic functions and \cite{MSRL} for some extensions to functions with singularities). In practice,   suitable approximation/compression techniques are to be employed.  Our preferred  option  is to look at the matrix $B$ as a rank-structured matrix  with the possibility  to   encode the structure by using numerical ranks.

Condensed representations  for rank-structured matrices have been proposed in a variety of papers.
Quasiseparable matrices are introduced in \cite{EG11}.  A complete review of their properties is presented in \cite{EGH_book}.  A quasiseparable representation of a matrix $A$ is defined  by two families of lower and upper generators  that are computed by exploiting the low-rank properties of the submatrices of $A$ located  in  its  lower and upper triangular part, respectively.  Given  in  input  a quasiseparable representation of $X_0$ and $B$  then the  structured  adaptation of \eqref{Newton} amounts to compute at each iteration  one  sum and  two products  of  quasiseparable matrices possibly
complemented with a  compression/approximation technique used to take under control the growth of the quasiseparable generators.   Generator-based algorithms   to perform these operations are described in
Chapter 4, 5 and 17 of \cite{EGH_book}.

A more flexible format  for rank-structured matrices   which is amenable for divide-and-conquer oriented techniques is called hierarchically semiseparable (HSS)  representation \cite{XCGL}.  This representation  is found by combining recursive partitioning, compression of off-diagonal blocks and  nestedeness  for the  generators of these blocks. In order to operate with HSS  matrices efficiently, one exploits
their representation with generators, demonstrated by the following
$4\times 4$ example:
\[
A=\left[\begin{array}{cccc}
   D_1 & G_1H_2 & G_1 R_1R_3 Q_3 & G_1 R_1 R_4 H_4\\
   P_2 Q_1 & D_2 & P_2 R_2 R_3 Q_3 & P_2 R_2 R_4 H_4\\
   G_3L_1 L_3 Q_1& G_3 L_1 L_4 H_2& D_3 & G_3 H_4 \\
   P_4 L_2 L_3 Q_1& P_4 L_2 L_4 H_2& P_4 Q_3 & D_4
\end{array}\right], 
\]
 where $D_i$ are square matrices of equal size. The representation is condensed if all the matrices $L_i$ and $R_i$ have sizes  less than a small constant $k\ll N$.  The generators $G_i$, $P_i$, $H_i$ and $Q_i$ are tall or skinny  matrices.  The value of $k$ is related with the maximum rank of all off-diagonal blocks at all levels of  the HSS recursive splitting of $A$ \cite{XCGL}.  Arithmetic operations   between two 
matrices of order $N$  expressed in a condensed HSS format can be performed in linear time w.r.t. $N$
\cite{XCGL}. 

\begin{figure}[ht] 
  \begin{subfigure}[b]{0.5\linewidth}
    \centering
    \includegraphics[width=0.75\linewidth]{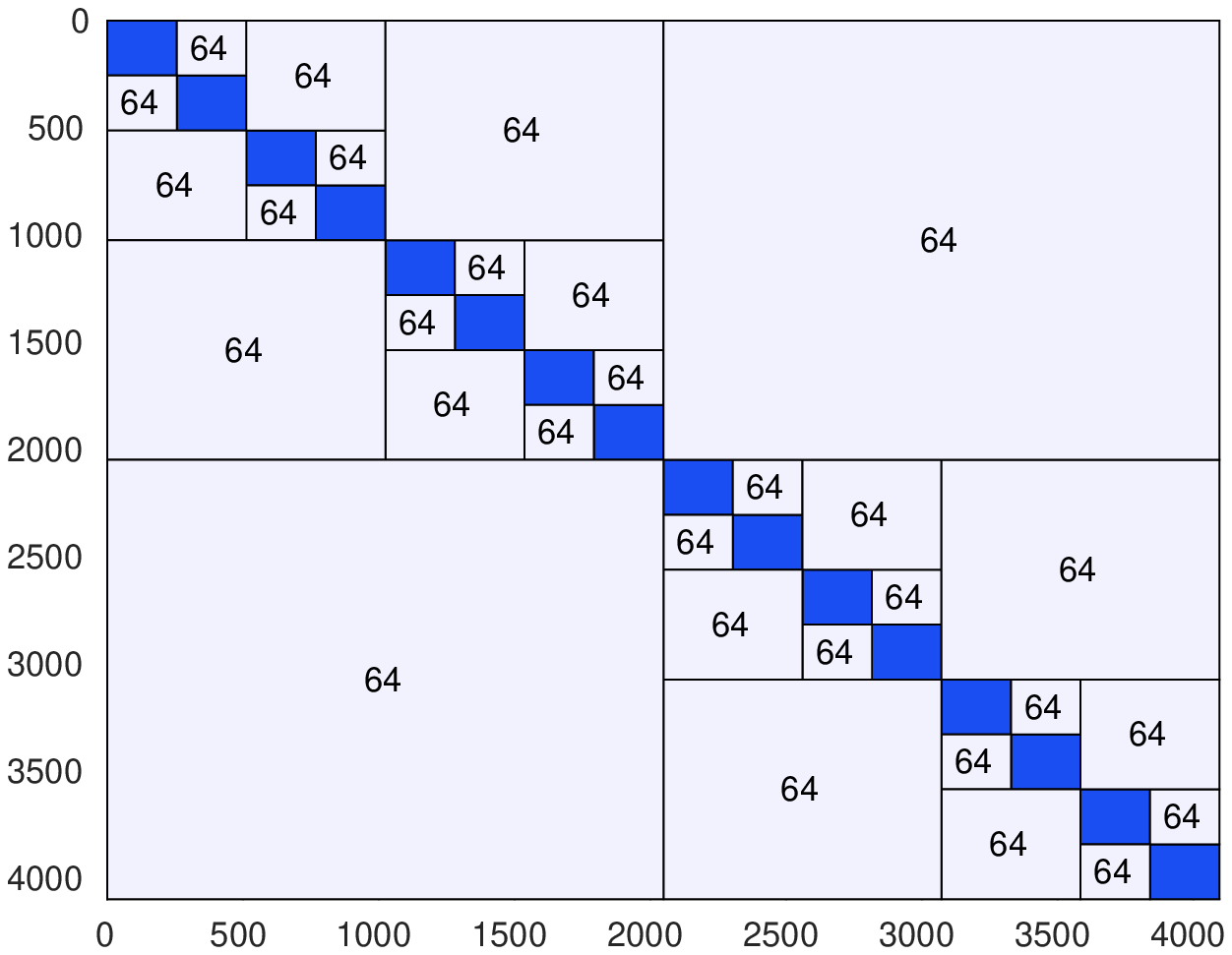} 
    \caption{Ranks of $A$ } 
    \label{fig7:a} 
    \vspace{4ex}
  \end{subfigure}
  \begin{subfigure}[b]{0.5\linewidth}
    \centering
    \includegraphics[width=0.75\linewidth]{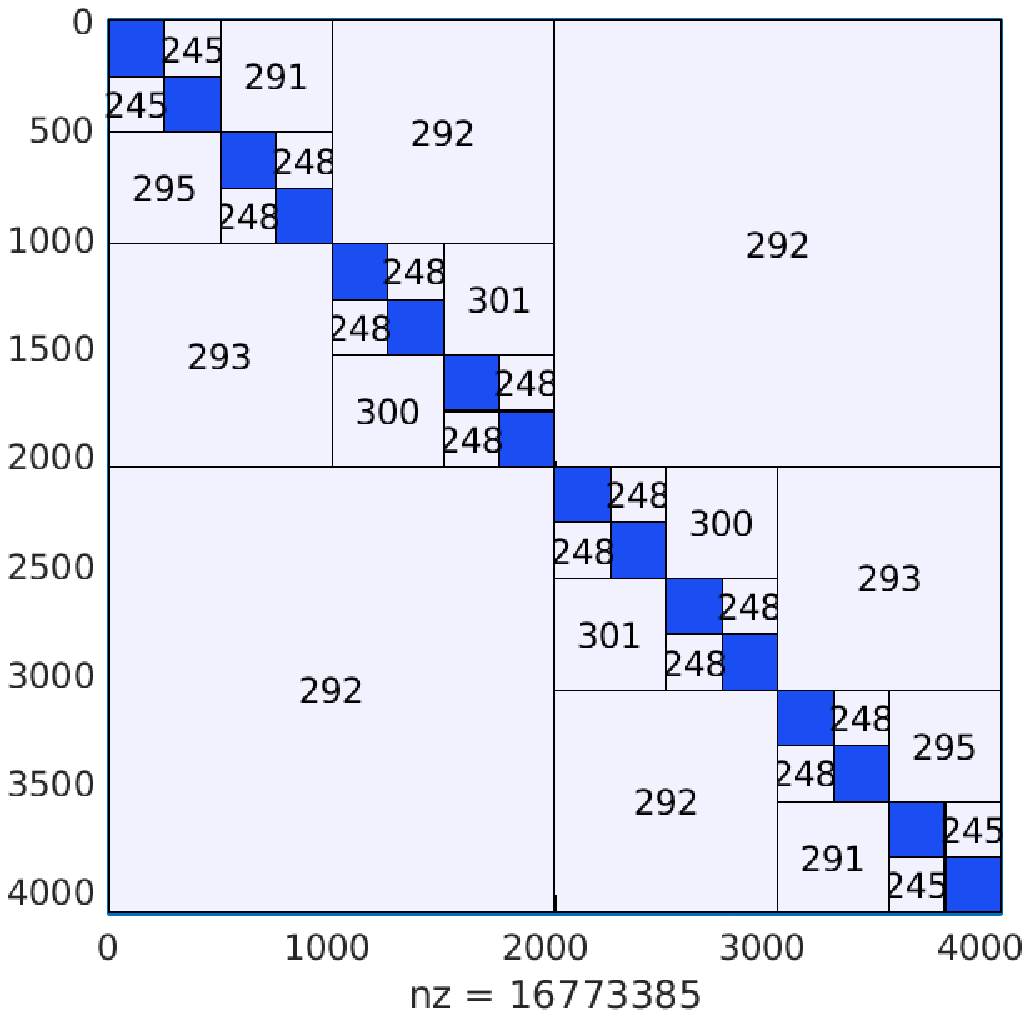} 
    \caption{Ranks of $X_0=(\phi_1(A))^{-1}$ } 
    \label{fig7:b} 
    \vspace{4ex}
  \end{subfigure} 
  \begin{subfigure}[b]{0.5\linewidth}
    \centering
    \includegraphics[width=0.75\linewidth]{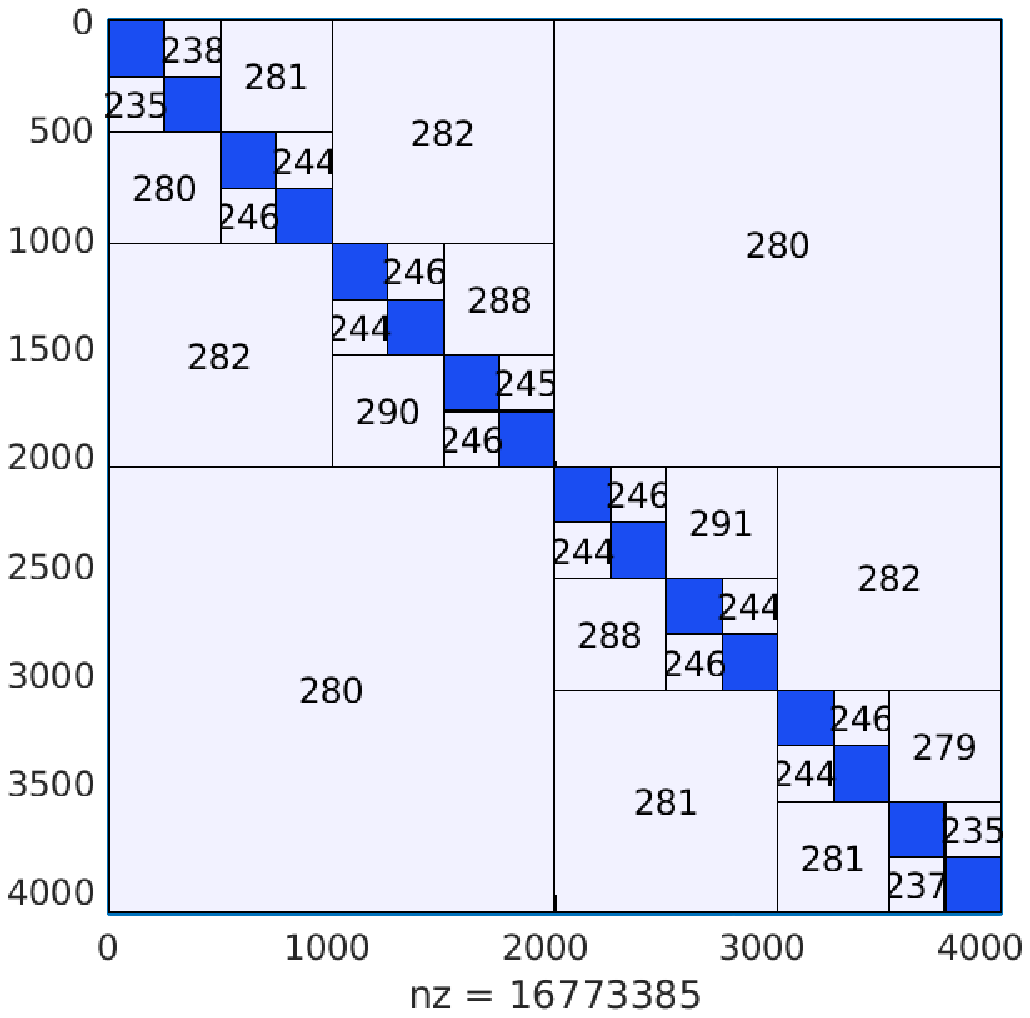} 
      \caption{Ranks of $\phi_2(A)$} 
    \label{fig7:c} 
  \end{subfigure}
  \begin{subfigure}[b]{0.5\linewidth}
    \centering
    \includegraphics[width=0.75\linewidth]{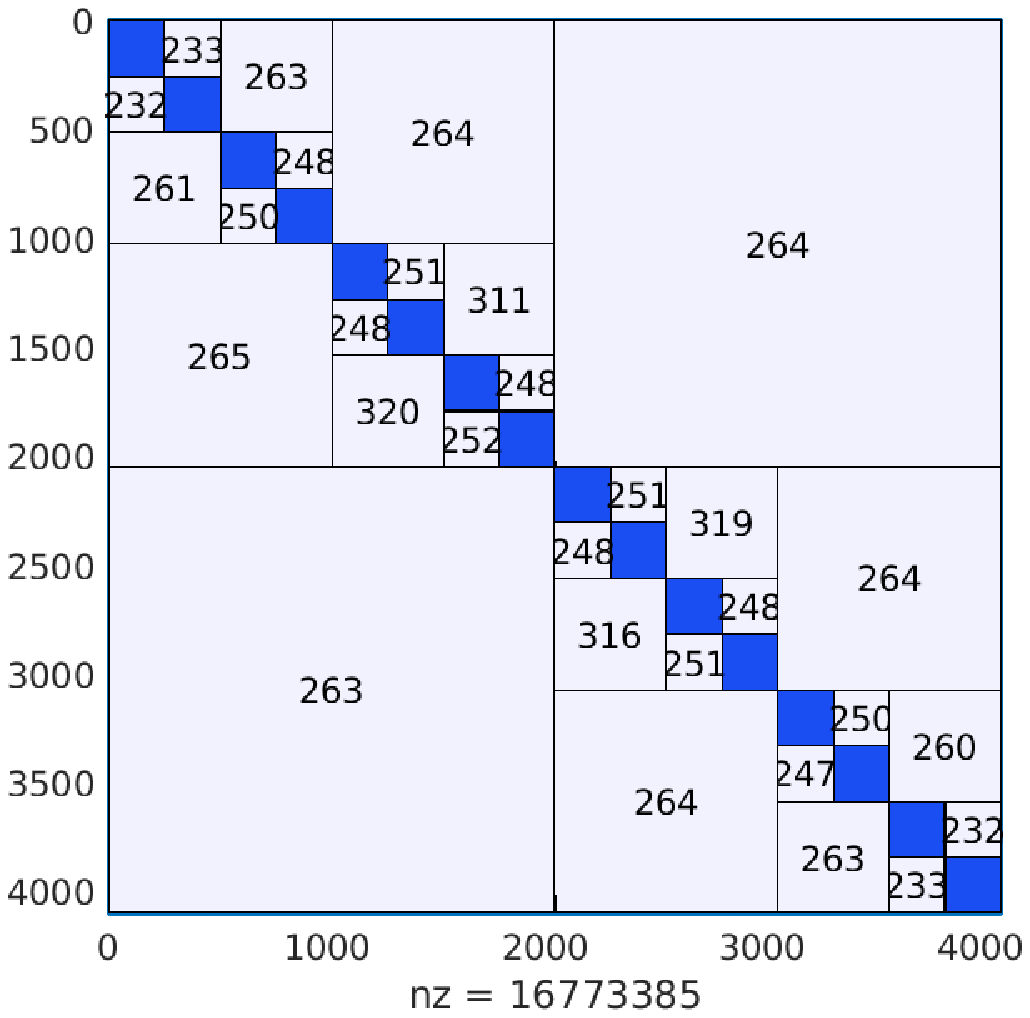} 
    \caption{Ranks of $X_7\simeq (\phi_2(A))^{-1}$ } 
    \label{fig7:d} 
  \end{subfigure} 
  \caption{Illustration of   the rank  properties of  HSS representations of the matrices involved in the computation of $(\phi_2(A))^{-1}$  by Newton's iteration  for the matrix $A$ given in \eqref{am}.}
  \label{fig7} 
\end{figure}

In view of the relation with the ranks of the off-diagonal blocks  it is clear that any arithmetic operation (except inversion) performed on HSS matrices can  increase their ranks.  Therefore, to use HSS   structure  efficiently under some iterative process we need  some compression algorithm. A Matlab toolbox to carry out  arithmetic operations among HSS matrices in exact or approximate  compressed form   is  described in \cite{MRK}.    Using this package, for the sake of illustration we show  in Figure  \ref{fig7} the   rank properties of the matrices generated by Newton's  iteration  applied for the computation of $(\phi_2(A))^{-1}$   starting from  $(\phi_1(A))^{-1}$     for a given rank-structured matrix $A$. Specifically, in our test we consider the matrix $A\in \mathbb R^{4096\times 4096}$   defined as follows:
\begin{equation}\label{am}
A=\frac{1}{3}\left[\begin{array}{cccc} M & N \\ N & \ddots  & \ddots \\ & \ddots & \ddots & N\\ & & N & M\end{array}\right]
\end{equation}
with 
\[
M=\left[\begin{array}{cccc} -8 & 1 \\ 1 & \ddots  & \ddots \\ & \ddots & \ddots & 1\\ & & 1 & -8\end{array}\right], \quad  N=\left[\begin{array}{cccc} 1& 1 \\ 1 & \ddots  & \ddots \\ & \ddots & \ddots & 1\\ & & 1 & 1\end{array}\right],
\]
and $M$ and $N$ of size $64$.
The matrix is generated in the solution of 2D Laplace's equation with Dirichlet boundary conditions by Q1 finite elements \cite{gondzio2011multilevel}.  In  Figure \ref{fig7} we show the numerical  ranks of the off-diagonal blocks  in the HSS representations of $A$, $(\phi_1(A))^{-1}$, $\phi_2(A)$ and  the approximation  $X_7$ of $(\phi_2(A))^{-1}$  generated by  Newton's iteration applied for the inversion
  of $\phi_{2}(A)$ with starting point $X_0=\psi_1(A)=(\phi_1(A))^{-1}$  stopped after $7$ iterations
  with error $\parallel \phi_{2}(A) X_7 -I \parallel_2\leq 2.3e-12$.  The compression threshold value is set to $1.0e-\!12$.

\section{A Krylov-type  Method for Computing the Action of  $\psi$-Functions on a  Vector}\label{threenew}           

The above results indicate  the possibility of approximating $B^{-1}=\psi_{\ell+1}(A)$, $\ell>0$,  using  Newton's method with starting point $X_0=\psi_{\ell}(A)$ provided that the eigenvalues of $A$ are properly localized. It follows easily from the arguments in Section \ref{two} that,   under the same assumption,  we can  apply  some Krylov-type  method  like GMRES  \cite{saadb} for approximating $\psi_{\ell+1}(A)\B b$  or, equivalently,
for solving the linear system \eqref{add1}. 
In this case $X_0$ plays the role of a preconditioner suitably determined to ensure the convergence of the projection method. 
The next result  immediately follows from Proposition 6.32 in \cite{saadb}.

\begin{proposition}\label{gmresp}
  Let $A\in \mathbb C^{N\times N}$   be a  diagonalizable matrix, i.e., $A=SDS^{-1}$, $D=\diag\left[\lambda_1, \ldots, \lambda_n\right ]$, with   all   eigenvalues $\lambda_i$, $1\leq i\leq N$,
  lying in the strip $ \Omega=\mathbb R\times \iu [-\pi/2,\pi/2]$ in the complex plane.  Let $\B x_m$, $m\geq 0$,  be the approximate solution of \eqref{add1},  with $B=\phi_{\ell+1}(A)$ and
  $X_0=\psi_{\ell}(A)$,  obtained from the $m$-th step of the GMRES algorithm, and let $\B r_m=X_0 \B b-X_0 B \B x_m$. Then, we have
  \[
  \parallel  \B r_m\parallel_2 \leq  \kappa_2(S) \rho(R)^m \parallel  \B r_0\parallel_2, \ \ R=I-X_0 B, \quad m\geq 0,
  \]
  where $ \rho(R)<1$  denotes the spectral radius of $R$ and $\kappa_2(S)$ is the 2-norm condition number  of $S$.
  \end{proposition}

  The convergence estimates for GMRES can be extended to general matrices by replacing the spectral decomposition of $A$ with its Jordan canonical form \cite{Sacchi}.   The bound in Proposition \ref{gmresp} depends on the 2-norm condition number of the eigenvector matrix $S$. This is satisfactory in the normal case but if  $A$ is far
from normal, then $\kappa_2(S)$ may have large magnitude and this    utterly invalidates  the bound.   Alternative GMRES convergence bounds based on the numerical range or the pseudospectra of the coefficient matrix have been proposed in the literature (see  \cite{EM}  and the references given therein).  However,  these bounds are not easy to compute and they  present other drawbacks so  that a common approach is to mitigate the impact of $\kappa_2(S)$ by assuming  that  the possible ill-conditioning is  due only to a low-dimensional invariant subspace  which contribution  can be deflated  in same way \cite{EM,SS}. 
  
  The computational cost
  for the GMRES algorithm  is dominated by the cost of  matrix-vector  multiplications  with the matrix $X_0 B$.  It is worth noticing that  from relation \eqref{recurrence}
  \begin{equation}\label{initmod}
  X_0 B =\psi_{\ell}(A)\phi_{\ell+1}(A)=A^{-1} (I-\psi_{\ell}(A)/\ell !), 
  \end{equation}
  which implies that the multiplication of $X_0 B $ by a vector   reduces to   first  multiply
  $\psi_{\ell}(A)$ by the same vector,  and then  solve a linear system with coefficient matrix $A$.   This is particularly interesting for $\ell=1$  since an efficient algorithm   to  evaluate $\psi_{1}(A)\B b$ has  been proposed in  \cite{BEG1, BEG2}.  The algorithm  relies upon  the  family of  polynomial/rational  expansions   of $\psi_{1}(z)$  given in \eqref{ratapps}. Based on \eqref{initmod},  complementing  the GMRES iterative solver with the  approximation  \eqref{ratapps}  provides an effective   method  for 
  computing the action of $\psi_{2}(A)$ on a vector.    A basic MatLab skeleton implementation is as follows: 
  
  \begin{algorithm}\label{a1} 
    \caption{\!\!\!: Given in input the matrix $A\in \mathbb C^{N\times N}$ and the vector $\B b\in \mathbb C^N$, this algorithm approximates the vector $\psi_{2}(A)\B b$}
    \label{algorithm2}
    \begin{algorithmic}[1] 
        \State \textbf{Select} the values of $n$ and $m$  in \eqref{ratappsnm}; 
        \State \textbf{Define} ${\tt funmv=@}(\B z)A\backslash(I_N-r_{n,m}(A))\B z$; 
        \State  \textbf{Call} $\B w= {\tt gmres}({\tt @funmv}, \B b, tol, maxit)$; 
            \State \textbf{return} $\B w$
    \end{algorithmic}
\end{algorithm}
  
 In view of  \eqref{initmod}, Algorithm 1 can be the building block  of a recursive procedure for the evaluation of $\psi_{\ell}(A)\B b$, $\ell>1$. Some explanations are, however,  in order with respect to the  execution of the first two steps. 
 
\begin{figure}[ht] 
  \begin{subfigure}[b]{0.5\linewidth}
    \centering
    \includegraphics[width=0.75\linewidth]{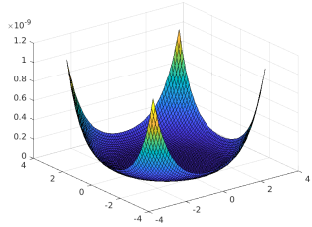} 
    \caption{$n=2$ and $m=32$ over $[-3,3]\times \iu[-3,3]$} 
    \label{fig11:a} 
    \vspace{4ex}
  \end{subfigure}
  \begin{subfigure}[b]{0.5\linewidth}
    \centering
    \includegraphics[width=0.75\linewidth]{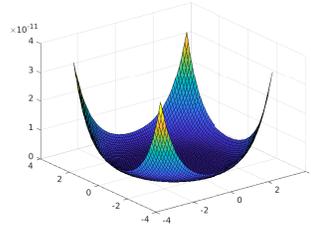} 
    \caption{ $n=2$ and $m=64$ over $[-3,3]\times \iu[-3,3]$}
    \label{fig11:b} 
    \vspace{4ex}
  \end{subfigure} 
  \begin{subfigure}[b]{0.5\linewidth}
    \centering
    \includegraphics[width=0.75\linewidth]{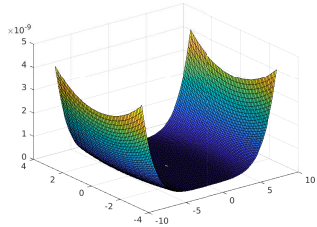} 
      \caption{$n=2$ and $m=64$ over $[-9,9]\times \iu[-3,3]$} 
    \label{fig11:c} 
  \end{subfigure}
  \begin{subfigure}[b]{0.5\linewidth}
    \centering
    \includegraphics[width=0.75\linewidth]{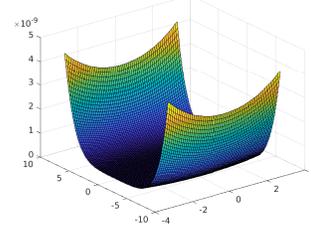} 
    \caption{$n=2$ and $m=64$ over $[-3,3]\times \iu[-9,9]$ } 
    \label{fig11:d} 
  \end{subfigure} 
  \caption{Surf plots of the  absolute errors $|\psi_{1}(z) -r_{n,m}(z)| $
  generated for different values of of $n$ and $m$ over different domains.}
  \label{fig11} 
\end{figure}

 F or a diagonalizable matrix $A$, the selection of $n$ and $m$ in \eqref{ratappsnm} is an approximation problem  depending on the localization of the spectrum of $A$.  In Figure \ref{fig11} we show  the plot of the absolute error $|\psi_{1}(z) -r_{n,m}(z)| $ for different values of $n$ and $m$ and  for different domains. 
 
 Recall that  $\psi_1(z)$ is a meromorphic function with poles $\pm 2 \pi\iu k$, $k\in \mathbb N$,  and a  removable singularity at  the origin.   It is remarkable that the approximation is quite accurate  even close around the singular points, whereas  the error increases appreciably with the size of the domain. A widespread approach to the computation of exponential and $\phi_\ell-$functions combines polynomial or Pad\'e approximation with a few steps of scaling-and-squaring \cite{higham2009scaling}. In principle, scaling-and-squaring may also be applied to our mixed polynomial-rational approximation, scaling the function argument by a suitable power of $2$ and then making use of the squaring formulas
\begin{equation}
\psi_1(2z)=\frac{2\psi_1(z)}{e^z+1}=\frac{2\psi_1(z)^2}{z+2\psi_1(z)}.\label{psisquaring}
\end{equation}
 The efficient implementation of   \eqref{psisquaring} is an ongoing research project and this  scaling-and-squaring scheme   is not considered here. 
 
 The  evaluation of the function ${\tt funmv}$ at step 2  of Algorithm 1 basically amounts to compute $r_{n,m}(A)\B z$. In the typical situation where $n\ll m$ this  computations reduces  to solve $m$  shifted linear systems of the form 
 \begin{equation}\label{lsyst}
 (A^2+(2 \pi k )^2 I_N)\B z_i=(A+2\pi \iu k I_N)(A-2\pi \iu k I_N)\B z_i=\B b, \quad 1\leq i\leq m.
 \end{equation}
 There is an extensive literature on the solution of shifted linear systems. In the case of interest where $A$ is quasiseparable we make use 
 of the  backward stable algorithm proposed in \cite{BEG1}. This algorithm  saves about half of computations in the solution of the  shifted linear systems by  reusing pieces of the structured QR factorization of the matrix $A$. 
 According to \cite{BEG1} for a quasiseparable matrix $A$ of size $N$ partitioned in blocks of size $n$ that are represented via quasiseparable generators of length $r\ll m$ the arithmetic  cost of solving the systems \eqref{lsyst} is of the order $4 n^2 m N$. 
 
 Some numerical tests showing the effectiveness  of  Algorithm 1    are presented in  Section \ref{five}.  
 
 \section{Numerical Results}\label{five}

We  have tested the application of Algorithm 1 for computing $\B w=\psi_{2}(A)\B b$ numerically by using MatLab. 

\begin{figure}[ht] 
  \begin{subfigure}[b]{0.5\linewidth}
    \centering
    \includegraphics[width=0.75\linewidth]{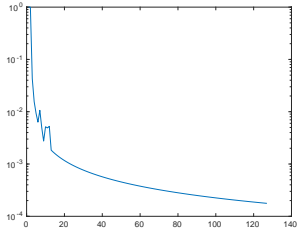} 
    \caption{Semilog plot of $err_1^{(j)}$ for the matrix in \ref{itnew1}
    } 
    \label{fig14:a} 
    \vspace{4ex}
  \end{subfigure}
  \begin{subfigure}[b]{0.5\linewidth}
    \centering
    \includegraphics[width=0.75\linewidth]{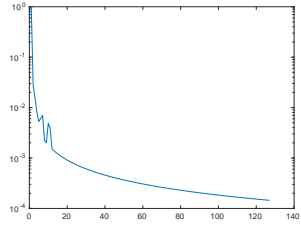} 
    \caption{Semilog plot of $err_2^{(j)}$ for the matrix in \ref{itnew1} }
    \label{fig14:b} 
    \vspace{4ex}
  \end{subfigure} 
  \begin{subfigure}[b]{0.5\linewidth}
    \centering
    \includegraphics[width=0.75\linewidth]{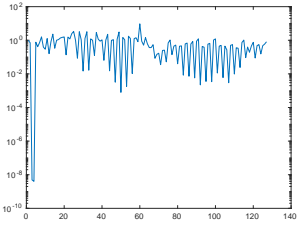} 
      \caption{Semilog plot of $err_1^{(j)}$ for the matrix in \ref{itnew2}} 
    \label{fig14:c} 
  \end{subfigure}
  \begin{subfigure}[b]{0.5\linewidth}
    \centering
    \includegraphics[width=0.75\linewidth]{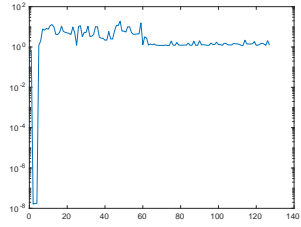} 
    \caption{Semilog plot of $err_2^{(j)}$ for the matrix in \ref{itnew2}} 
    \label{fig14:d} 
  \end{subfigure} 
  \caption{Semilog plots of  $err_1^{(j)}$ and $err_2^{(j)}$ for ther matrices  defined in \ref{itnew1} and \ref{itnew2} of size $N=128$ with $\epsilon=1.0e\!-14, 1.0e\!-8$, respectively.}
  \label{fig14} 
\end{figure}

Numerical experiments  have been  carried out  for comparison with the  classical approach 
based on the Arnoldi method \cite{DK,Kni,GS},  where $\B w_j=W_j \psi_{2}(H_j)\B e_1 
\parallel \B b\parallel_2$, $j\geq 1$,  is taken as an approximation of $\B w$ and $W_j$ and $H_j$ are generated in the  Arnoldi  process.  This scheme works quite well in general, notwithstanding that $\psi_2(z)$ is a meromorphic  function.  The crux is that  the performance is depending on a number of issues such as the choice of the starting vector and the stopping criterion in the Arnoldi process as well as the properties of the spectrum of the matrices $H_j$ and the quality of the polynomial approximation of $\psi_2(z)$ on this spectrum.  These issues  can be  difficult to tackle and  resolve for a  class of matrices.  To see this let us consider the following examples: 
\begin{enumerate}
    \item\label{itnew1} $A=Z+\epsilon \B e\B e^T$, $Z=\left[\begin{array}{cccc} & 1 \\ && \ddots \\ &&&1 \\1 \end{array}\right]$, $\B e^T=\left[1, \ldots, 1\right]$; 
    \item\label{itnew2}  $A=\diag(-1-{\tt lambertw}((-1)*[1,1]-1,-{\tt exp}(-1)))\otimes I + +\epsilon \B e\B e^T$, 
\end{enumerate}
where ${\tt lambertw}$ computes the Lambert W function \cite{CGH}  which is involved in the  numerical computation of the poles of $\psi_2(z)$ \cite{KT}. We have implemented the Arnoldi-based method in  MatLab. The approximation $\B w_j$ is computed using ${\tt expm}$ and the backslash operator.  As a stopping criterion we evaluate the relative error $err_1^{(j)}=\parallel \B w_{j+1} - \B w_j\parallel_2/\parallel \B w_j\parallel_2$. As a measure of accuracy we also  determine the  relative error $err_2^{(j)}=\parallel \B w- \B w_j\parallel_2/\parallel \B w\parallel_2$, where $\B w$ is computed  in some  way (varying with the considered test).  Specifically,  the matrix in the first test is well conditioned  with eigenvalues far from the poles of $\psi_2(z)$ and $\B w$ is found using ${\tt expm}$ and the backslash operator applied to $\phi_2(A) \B w=\B b$ with $\B b =\B e_1$. On the contrary,  the eigenvalues of the matrix in \ref{itnew2} are clustered around the first two poles of $\psi_2(z)$ so that $\B w$  is fixed equal to $\B e$ and then $\B b$ is determine by $\phi_2 \B w=\B b$.  In Figure \ref{fig14} we show the plots of $err_1^{(j)}$ and $err_2^{(j)}$.  The matrices  have size $N=128$ and the value of $\epsilon$ is set to $1.0e\!-14$ and $1.0e\!-8$ in the first and second example, respectively. 
We notice that the Arnoldi process applied to the matrix $A$ defined in \ref{itnew1} generates ill-conditioned submatrices $H_j$ that have many eigenvalues clustered around the origin in the complex plane.  The effect is a deterioration of the  precision.   Perhaps,  the effect   could be alleviated  replacing the backslash operator with Newton's iteration,  but this would need  some localization of the spectrum of the matrices  $H_j$. 
For the matrix in \ref{itnew2} we observe a similar behaviour of the spectrum of  the matrices $H_j$ which accumulates around the two poles and the origin.  In a few steps the errors reach a minimum value of order $1.0e\!-8$ and after that rapidly grow and  stabilize around the unit. 

For comparison in Table \ref{tnew1} and \ref{tnew2} we  describe the results obtained by Algorithm 1 applied to the matrix in \ref{itnew1}  and \ref{itnew2}.  For each example we show the value  of  $n$   and $m$ in the rational approximant, the number $it_{\tt gmres}$ of iterations  of {\tt gmres}, the relative residual $r_{\tt gmres}$ of the approximation returned by {\tt gmres}  together  with the relative error $err_2$. 

\noindent
\begin{table}
		\begin{center}
			\begin{tabular}{|l|l|l|l|} \hline
			
			$(n,m)$ & $(3,8)$ &$ (3,16)$ & $(3,32)$\\ 	\hline
			$it_{\tt gmres}$ & $17$ & $17$ & $17$  \\ \hline
			$r_{\tt gmres}$ &  $4.8e\!-14$ &  $ 4.8e\!-14$ & $ 4.8e\!-14$  \\ \hline
			 $err_2$&  $8.7e\!-14$ &  $ 5.3e\!-14$ & $ 5.3e\!-14$  \\
			\hline
			\end{tabular}
\end{center}
       \caption{Performance of Algorithm 1 applied to the matrix in \ref{itnew1}  of size $N=128$ with $\epsilon=1.0e\!-14$ and  the tolerance of {\tt gmres} set at $tol=1.0e\!-12$. \label{tnew1}.} 
\end{table}

     \noindent
\begin{table}
		\begin{center}
			\begin{tabular}{|l|l|l|l|l|} \hline
			
			$(n,m)$ & $(3,16)$ &$ (3,32)$ & $(3,64)$ & $(3,128)$\\ 	\hline
			$it_{\tt gmres}$ & $3$ & $3$ & $3$ & $3$  \\ \hline
			$r_{\tt gmres}$ &  $1.8e\!-15$ &  $ 1.8e\!-15$ & $ 1.7e\!-15$ &  $ 1.8e\!-15$\\ \hline
			 $err_2$&  $5.0e\!-3$ &  $ 4.4e\!-5$ & $3.6e\!-7$ &  $1.8e\!-8$ \\
			\hline
			\end{tabular}
\end{center}
       \caption{Performance of Algorithm 1 applied to the matrix in \ref{itnew2}  of size $N=128$ with $\epsilon=1.0e\!-8$ and  the tolerance of {\tt gmres} set at $tol=1.0e\!-12$. \label{tnew2}.} 
\end{table} 

\begin{figure}
  \centering
  \subfloat[]{\includegraphics[width=0.5\textwidth]{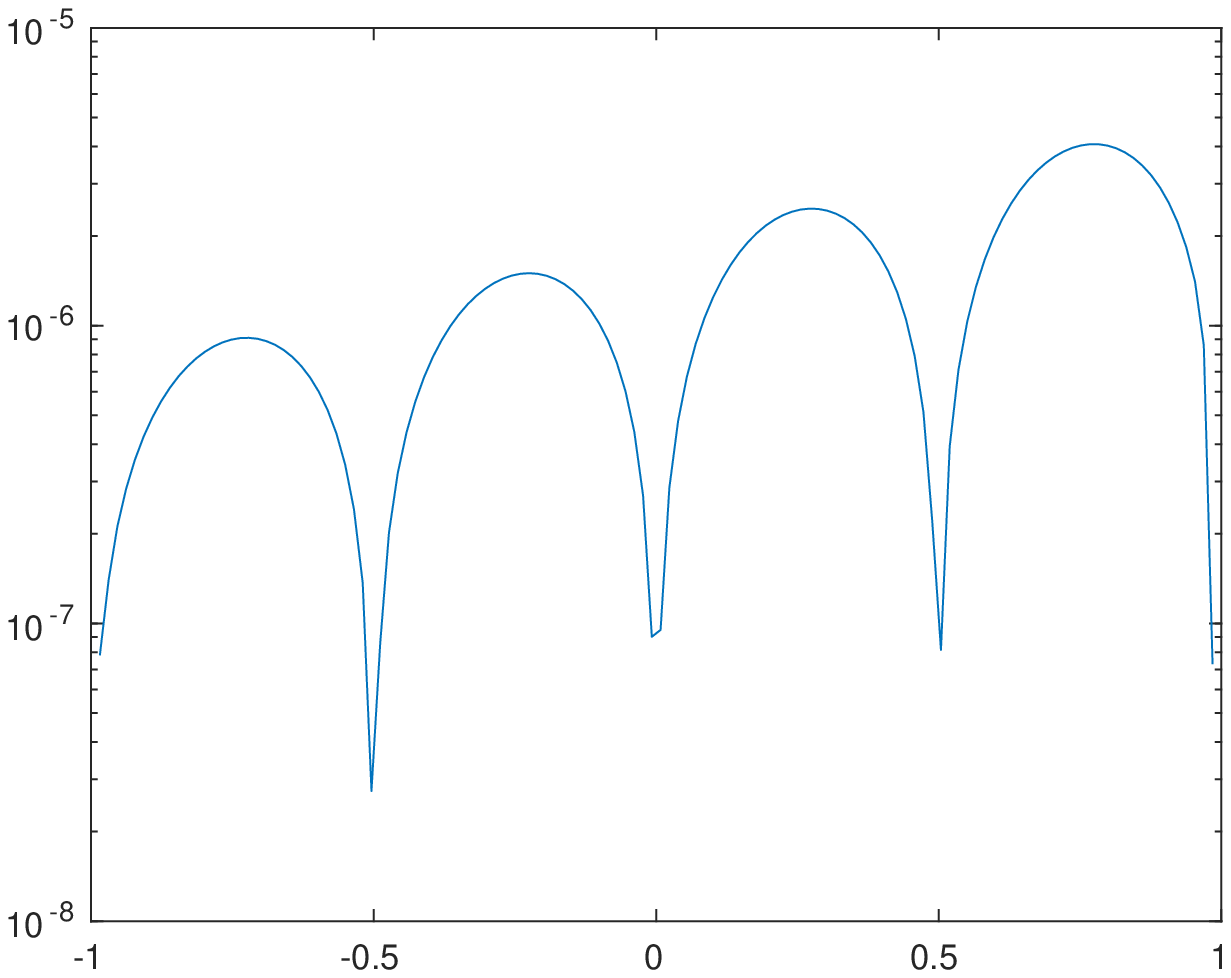}}
  \hfill
  \subfloat[]{\includegraphics[width=0.5\textwidth]{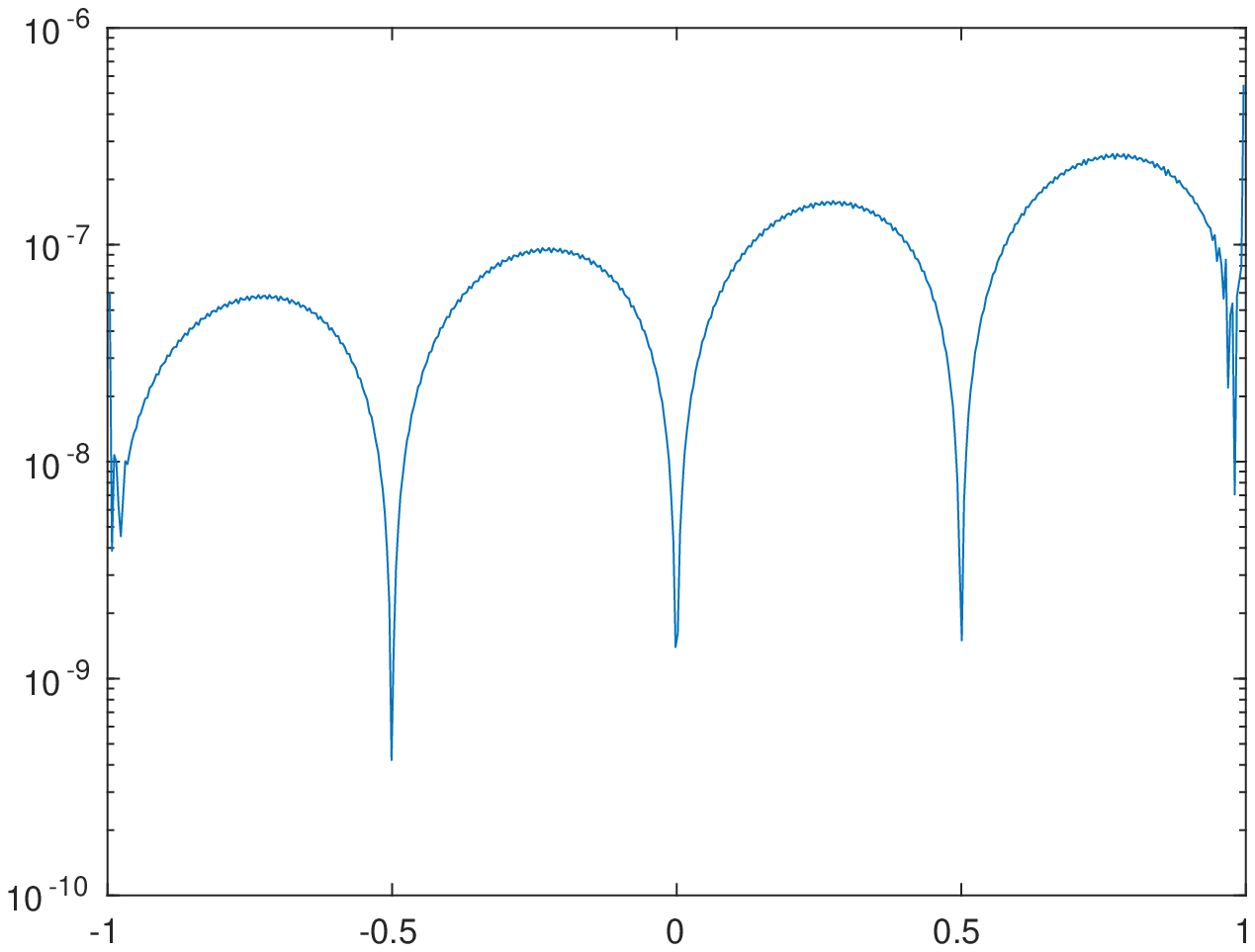}}
  \caption{Plots of the error vectors generated  by Algorithm 1  for  $N=128$ (a) and $N=512$  (b). }
\label{fferr}
\end{figure}

Numerical tests have been also performed to investigate the application of Algorithm 1 in the solution of  the inverse  problems  described in the introduction.  For the sake of illustration let us  consider the following differential problem: 
\begin{equation}\label{heat2}
\frac{\partial u(z,t)}{\partial t}=\frac{e^{z-4}}{\sigma^2}\frac{\partial^2 u(z,t)}{\partial z^2} + tf(z) , \quad f(z)= \sin(2\pi z), \ (z,t)\in [-1,1]\times [0,1],
\end{equation}
with boundary  conditions $u(-1, t)=u(1,t)=0$, $u(z, 0)=0$ and $\sigma=10$. The differential problem  has been  solved in Mathematica using the internal function NDSolve with extended precision.  The computed solution $u(z,t)$ evaluated at $t=1$  is used to define $h(z)=u(z,1)$.  Then the inverse problem concerns the reconstruction of $f(z)$ from the boundary conditions and the additional constraint  $u(z,1)=h(z)$. Using a discretization in space by finite differences over $N+2$ equispaced points in the interval $[-1,1]$ leads  to the first order system 
 \[
  \odv{\B u}{t}= A \B u(t)+ t \B {f},
  \]
  \[
  A=\left(\frac{N+1}{2 \sigma}\right)^2 \diag(e^{z_1-4}, \ldots, e^{z_N-4})\left[\begin{array}{cccc}-2 & 1\\ 1& \ddots & \ddots \\& \ddots & \ddots & 1\\ & & 1& -2\end{array}\right]
  \]
  with conditions
  \[
   \B u(0)=\B 0,  \quad  \B u(1)=\B h,  \]
   and  
$\B f=\left[ f(z_1), \ldots, f(z_N)\right]^T$, $\B h=\left[ h(z_1), \ldots, h(z_N)\right]^T$.  The unknown  vector $\B f$ can thus be determined  by means of formula \eqref{invex}. The matrix 
  $A$ is  similar  to a  negative definite  matrix  and therefore our methods can be applied. 
In  Figure \ref{fferr}  we  plot the absolute error vector with entries  $|\hat f_i -f(z_i)|$, $1\leq i\leq N$, $N\in \{128, 512\}$,   where  $\hat f_i$  are generated by  Algorithm 1 with $n=2$, $m=32$, $tol=1.0e\!-10$ and $maxit=40$. 
 The {\tt {gmres}} command  detects convergence at iteration  7 and 8 for  $N=128$ and $N=512$, respectively.  The  finer discretization produces a small error. Similar plots are observed for the vectors generated by  using {\tt {expm}}  and the backslash operator.  The condition number of the matrices involved is of order $1.0e\!+8$.

\section{Conclusions  and Future Work}\label{six}
In this paper we have presented two approaches  based on Newton's iteration and Krylov-type methods
for the  efficient computation of the inverse  of  a matrix $\phi$-function as well as the action of this inverse matrix  on a vector.  In particular,  an appealing iterative procedure for  computing   $\psi_2(A)\B v$ has been devised. Numerical experiments show  that the proposed methods  exhibit
good robustness and convergence  properties.  The  iterative scheme for the  approximation of  $\psi_2(A)\B v$  requires at each step  to compute an approximation of  products of the form $\psi_1(A)\B w$ by solving several linear systems whose
matrices differ from A by a complex  multiple of the identity matrix.  Future work is concerned with  the efficient solution of these  shifted systems using the techniques introduced in \cite{BPD}. Another  interesting research topic  would be the design of an adaptive modification of the algorithm in \cite{BEG1, BEG2} for computing $\psi_1(A)$ capable  to  determine automatically the "best"  polynomial/rational approximation formula for $\psi_1(z)$. 

\section*{Acknowledgment}
The author would like to  thank Prof. Paola Boito and Prof. Yuli Eidelman for useful discussions and feedback. The author is  also indebted to  Yuli Eidelman for the  English  translation of  reference \cite{KT}.


\end{document}